\documentclass{amsart}
\usepackage{amsmath}
\usepackage{amssymb}
\usepackage{amsfonts}
\usepackage{amsthm}
\usepackage{enumerate}
\usepackage[all]{xy}
\usepackage[latin1]{inputenc}
\usepackage{mathdots}
\usepackage{graphicx}
\usepackage{latexsym}
\usepackage{mathrsfs}
\usepackage{textcomp}
\usepackage{color}
\input xy

\makeatletter

\newtheoremstyle{definition}
        {5pt}
        {3pt}
        {}
        {0pt}
        {\scshape}
        {.}
        {5pt}
        {\thmname{#1} \thmnumber{#2} \thmnote{[#3]}} 

\newtheoremstyle{theorems}
        {5pt}
        {3pt}
        {\itshape}
        {0pt}
        {\scshape}
        {.}
        {5pt}
        {\thmname{#1} \thmnumber{#2}\thmnote{[#3]}} 

\swapnumbers

\theoremstyle{theorems}
\newtheorem{Theo}{Theorem}[section]
\newtheorem{TheoR}[Theo]{Theorem \cite{Rin}}

\newtheorem{Prop}[Theo]{Proposition}

\newtheorem{Lemma}[Theo]{Lemma}

\theoremstyle{definition}
\newtheorem{Defn}[Theo]{Definition}

\newtheorem{DefR}[Theo]{Definition \cite{Rin}}

\newcommand{\Sa}{\mathit\Sigma}
\newcommand{\Ga}{\mathit\Gamma}
\newcommand{\Da}{\mathit\Delta}

\newcommand{\GaA}{\Ga\hspace{-1.5pt}_A}

\newcommand{\Hom}{{\rm Hom}}
\newcommand{\Ext}{{\rm Ext}}
\newcommand{\End}{{\rm End}}

\newcommand{\mmod}{{\rm mod}\hspace{0.4pt}}
\newcommand{\ind}{{\rm ind}\hspace{0.4pt}}
\newcommand{\rad}{{\rm rad}}

\newcommand{\taum}{{\tau^-\hspace{-1pt}}}

\makeatother

\newcommand{\Z}{{\mathbb{Z}}}

\begin{document}

\title[Tilted algebras]{\sc Another characterization of tilted algebras}



\author[Shiping Liu]{Shiping Liu}

\address{Shiping Liu\\ D\'epartement de math\'ematiques, Universit\'e de Sherbrooke, Sherbrooke, Qu\'ebec, Canada.}
\email{shiping.liu@usherbrooke.ca}

\thanks{This research is supported in part by the Natural Sciences and Engineering Research Council of Canada.}

\maketitle

\begin{abstract}

\vspace{15pt}

We give a new characterization of tilted algebras by the existence of certain special subquivers in their Auslander-Reiten quiver. This result includes the existent characterizations of this kind and yields a way to obtain more tilted quotient algebras from a given algebra.

\end{abstract}

\bigskip
\medskip

\section*{Introduction}

\medskip

Since its introduction by Happel and Ringel, the theory of tilted algebras has been one of the most important topics in the representation theory of artin algebras; see, for example, \cite{HaR, Liu5, Ker, Rin, Rin3}. Indeed, this class of algebras is closely related to hereditary algebras, a relatively well-understood class of algebras. Initially, tilted algebras are characterized by the existence of a slice in their module category, or equivalently, the existence of a slice module; see \cite[(4.2)]{Rin}. Later, the convexity of a slice module was replaced by a weaker condition; see \cite{Bak}. Most recently, a slice module is replaced by a sincere module which is not the middle term of any short chain in the module category; see \cite{JMS}.

\medskip

All the above-mentioned characterizations of tilted algebras require some know\-ledge of the entire module category, and hence, they are rather difficult to be verified for algebras of infinite representation type. To overcome this difficulty, replacing the convexity of a slice with respect to arbitrary maps by the convexity with respect to irreducible maps, one characterizes tilted algebras by the existence of an Auslander-Reiten component which contains a faithful section admitting no backward maps to its Auslander-Reiten translate; see \cite{Liu3, Sko2}. Later, a section in this characterization is replaced by a slightly weaker notion of a left section; see \cite{Ass}.

\medskip

Observe that a section, as well as a left section, requires some knowledge of an entire Auslander-Reiten component. In this paper, by relaxing the convexity of a slice in the module category, we obtain a locally defined notion of a {\it cut} in the Auslander-Reiten quiver. Our main result says that an artin algebra is tilted if and only if its Auslander-Reiten quiver contains a faithful cut admitting no backward map to its Auslander-Reiten translate. This not only is easy to be verified but also yields a way to obtain more tilted quotient algebras from a given artin algebra.

\medskip

\section{Preliminaries}

\medskip

Throughout this paper, $A$ stands for an artin algebra. We shall denote by $\mmod A$ the category of finitely generated left $A$-modules, by $\ind A$ the full subcategory of $\mmod A$ generated by the indecomposable modules. Recall that a module $T$ in $\mmod A$ is called {\it tilting} \hspace{0.3pt} if ${\rm pdim}(T)\le 1$, $\Ext_A^1(T, T)=0$, and the number of non-isomorphic indecomposable direct summands of $T$ is equal to the number of non-isomorphic simple $A$-modules.
The {\it Jacobson radical} of $\mmod A$, written as ${\rm rad}(\mmod A)$, is the ideal generated by the non-invertible maps between indecomposable modules, while the {\it infinite radical} ${\rm rad}^\infty(\mmod A)$ is the intersection of all powers ${\rm rad}^n(\mmod A)$ with $n\ge 0$. A map in $\rad(\mmod A)$ will be called a {\it radical} map.

\medskip

We shall use freely some standard terminology and some basic results of the Auslander-Reiten theory of irreducible maps and almost split sequences, for which we refer to \cite{ARS}. As usual, we shall denote by $\GaA$ the Auslander-Reiten quiver of $\mmod A$. This is a translation quiver whose vertex set is a complete set of representatives of the isomorphism classes of modules in $\ind A$, and whose arrows correspond to irreducible maps, and whose translation is given by the Auslander-Reiten translations $\tau={\rm DTr}$ and $\tau^-={\rm TrD}$, where $D$ denotes the standard duality for $\mmod A$. For simplicity, we shall write $\tau X=0$ if $X$ is projective and $\taum X=0$ if $X$ is injective.

\medskip

Let $\Sa$ be a full subquiver of $\GaA$. The {\it annihilator} of $\Sa$, written as ${\rm ann}(\Sa)$, is the intersection of all annihilators ${\rm ann}(M)$ with $M\in \Sa$. One says that $\Sa$ is {\it faithful} if ${\rm ann}(\Sa)=0$ and {\it sincere} if every simple $A$-module is a composition factor of some module in $\Sa$. Recall that $\Sa$ is a {\it section} in a connected component $\Ga$ of $\GaA$ if $\Sa$ is a connected subquiver of $\Ga$, which
contains no oriented cycle, meets each $\tau$-obit in $\Ga$ exactly once, and is convex in $\Ga$, that is, every path in $\Ga$ with end-points belonging to $\Sa$ lies entirely in $\Sa$; see \cite[(2.1)]{Liu1}.

\medskip

Recall that a {\it path} in $\ind A$ is a sequence of non-zero radical maps
$$\xymatrixcolsep{22pt}\xymatrix{X_0\ar[r]^{f_1} & X_1 \ar[r] & \cdots \ar[r] & X_{n-1}\ar[r]^{f_n}& X_n}
\vspace{1pt}$$
between modules in $\ind A$; and such a path is called {\it non-zero} if $f_1\cdots f_n\ne 0.$ The following result is implicitly included in the proof of \cite[(1.2)]{Liu3}.

\medskip

\begin{Lemma}\label{inf-rad} Let $X, Y$ be modules in $\ind A$ such that $\rad^\infty(X, Y)\ne 0.$ For each integer $n>0$, $\ind A$ contains a non-zero path
$$\xymatrixcolsep{20pt}\xymatrix{X\ar[r]^{u_n}& Y_n\ar[r]^{f_n} & Y_{n-1}\ar[r] & \cdots \ar[r] & Y_1 \ar[r]^{f_1} & Y,}
\vspace{-1pt}$$ where $u_n\in \rad^\infty(X, Y_n)$ and $f_n, \cdots, f_1$ are irreducible$\hspace{0.5pt};$ and a non-zero path
$$
\xymatrixcolsep{20pt}\xymatrix{X \ar[r]^-{g_1} & X_1\ar[r] & \cdots \ar[r] & Y_{n-1}\ar[r]^-{g_n} &  X_n \ar[r]^{v_n}&Y,}
\vspace{-1pt}$$ where $g_1, \cdots, g_n$ are irreducible and $v_n\in \rad^\infty(X_n, Y)$.

\end{Lemma}

\medskip

\begin{Defn}\label{map-depth}

Let $f: X\to Y$ be a map in ${\rm mod}\hspace{0.5pt}A$. We define the {\it depth} of $f$, written as ${\rm dp}(f)$, to be infinity in case $f\in {\rm rad}^\infty(X, Y)$; and otherwise, to be the integer $n\ge 0$ for which $f\in {\rm rad}^n(X, Y)$ but $f\notin {\rm rad}^{n+1}(X, Y)$.

\end{Defn}

\medskip

\noindent{\sc Remark.} A map $f: X\to Y$ in $\ind A$ is irreducible if and only if ${\rm dp}(f)=1$.

\medskip

The following result is well known.

\medskip

\begin{Lemma}\label{fdp-map}

Let $f: X\to Y$ be a radical map in $\ind A$. If ${\rm dp}(f)<\infty$, then $$f=f_1+\cdots +f_n+g,\vspace{-2pt}$$
where $g\in \rad^\infty(X, Y)$, and $f_1, \ldots, f_n$ are non-zero composites of irreducible maps between modules in $\ind A$.

\end{Lemma}

\medskip

A path $\xymatrixcolsep{18pt}\xymatrix{X_0\ar[r] & X_1 \ar[r] & \cdots \ar[r] & X_n}$ in $\GaA$ is {\it sectional} if $X_{i-1}\ne \tau X_{i+1}$ for every $0<i<n;$ {\it presectional} if there exists an irreducible map $u_i: \tau X_{i+1}\oplus X_{i-1}\to X_i$ for each $0<i<n;$ and {\it non-zero} if there exist irreducible maps $f_i: X_{i-1} \to X_i$, $i=1, \ldots, n,$ such that $f_1\cdots f_n\ne 0$. It is well known that a sectional path is presectional. For convenience, we quote the following result from \cite[(1.15)]{Liu1}.

\medskip

\begin{Lemma}\label{pre-sec-path}

If $\xymatrixcolsep{16pt}\xymatrix{X_0\ar[r] & X_1 \ar[r] & \cdots \ar[r]& X_n}$ is a presectional path in $\GaA,$ then there exist irreducible maps $f_i: X_{i-1} \to X_i,$ $i=1, \ldots, n,$ such that ${\rm dp}\hspace{0.4pt}(f_n\cdots f_1)=n.$

\end{Lemma}

\medskip

Applying first the Harada-Sai Lemma; see \cite{HS} and then Lemma \ref{pre-sec-path}, we obtain immediately the following result.

\medskip

\begin{Lemma}\label{lem-psp}

If $\Sa$ is a finite subquiver of $\Ga_A$, then the lengths of the non-zero paths in $\ind A$ passing through only modules in $\Sa$ are bounded, and consequently, $\Sa$ has no infinite sectional path.

\end{Lemma}

\medskip

The weak convexity defined below is essential for our investigation.

\medskip

\begin{Defn}

Let $\Sa$ be a full subquiver of $\GaA$. We shall say that $\Sa$ is

\vspace{-1.5pt}

\begin{enumerate}[$(1)$]

\item {\it convex} in $\ind A$ if every path in $\ind A$ with end-points in $\Sa$ passes through only modules in $\Sa;$ and

\item {\it weakly convex} in $\ind A$ if every non-zero path in $\ind A$ with end-points in $\Sa$ passes through only modules in $\Sa.$

\end{enumerate}

\end{Defn}

\medskip

\begin{Lemma}\label{lem-wcm}

Let $\Sa$ be a finite subquiver of $\GaA$, which is weakly convex in $\ind A$.

\begin{enumerate}[$(1)$]

\item If $X, Y\in \Sa$, then $\rad^\infty(X, Y)=0.$

\vspace{1pt}

\item The endomorphism algebra of the direct sum of the modules in $\Sa$ is connected if and only if $\Sa$ is connected.

\end{enumerate}

\end{Lemma}

\noindent{\it Proof.} (1) Suppose that there exists a non-zero map $f\in \rad^\infty(X, Y),$ for some $X, Y\in \Sa$. For each integer $n>0$, by Lemma \ref{inf-rad}, $\ind A$ has a non-zero path
$$\xymatrixcolsep{20pt}\xymatrix{X\ar[r]^{g_n} & Y_n\ar[r]^-{f_n} & Y_{n-1}\ar[r] & \cdots \ar[r] & Y_1 \ar[r]^{f_1} & Y,}$$ where the $f_i$ are irreducible. Since $\Sa$ is weakly convex, $Y_1, \ldots, Y_n\in \Sa$, a contradiction to Lemma \ref{lem-psp}.

(2) Suppose that the endomorphism algebra of the direct sum of the modules in $\Sa$ is connected. Let $X, Y\in \Sa$ be distinct modules. Then $\Sa$ contains modules $X=X_1, X_2, \ldots, X_n=Y$ such that, for each $1\le i<n$, there exists a non-zero map $f_i$ from $X_i$ to $X_{i+1}$ or from $X_{i+1}$ to $X_i$. By Statement (1), ${\rm dp}(f_i)<\infty;$ and by Lemma \ref{fdp-map}, $\GaA$ has a non-zero path $\rho_i$ from $X_i$ to $X_{i+1}$ or from $X_{i+1}$ to $X_i$, for all $1\le i<n$. Since $\Sa$ is weakly convex, all the paths $\rho_i$ with $i1\le i<n$ lie in $\Sa$. This shows that $\Sa$ is connected. The proof of the lemma is completed.

\medskip

Recall that $A$ is {\it tilted} if $A=\End_H(T)$, where $H$ is a hereditary artin algebra and $T$ is a tilting module in $\mmod H$; see \cite{HaR}. It is a well-known result of Ringel's that $A$ is tilted if and only if $\mmod A$ contains a slice; see \cite[(4.2)]{Rin}. Observe that a slice in $\mmod A$ is precisely the full additive subcategory of $\mmod A$ generated by the modules isomorphic to those in a {\it slice} of $\GaA$ as defined below.

\medskip

\begin{DefR}\label{def-slice}

A full subquiver $\Da$ of $\GaA$ is called a {\it slice} if it satisfies the following conditions.

\vspace{-1pt}

\begin{enumerate}[(1)]

\item The subquiver $\Da$ is sincere and convex in $\ind A$.

\item If $X\in \Da$, then $\tau X\notin \Da$.

\item If $X\to Y$ is an arrow in $\GaA$ with $Y\in \Da$, then either $X$ or $\taum X$ belongs to $\Da$.

\end{enumerate} \end{DefR}

\medskip

For convenience, we reformulate Ringel's result as follows. Although it is stated  in \cite[(4.2)]{Rin} for a finite dimensional algebra over an algebraically closed field, the same proof work for an artin algebra.

\medskip

\begin{TheoR}\label{Lem-slice}

Let $A$ be an artin algebra, and let $\Da$ be a full subquiver of $\GaA$.

\vspace{-1pt}

\begin{enumerate}[$(1)$]

\item If $A=\End_H(T)$ with $H$ hereditary and $T$ a tilting $H$-module, then $T$ determines a slice in $\GaA$ generated by the direct summands of $\Hom_A(T, D(H))$.

\vspace{1pt}

\item The subquiver $\Da$ is a slice if and only if $S=\oplus_{X\in \Da}X$ is a tilting module in $\mmod A$ such that $H=\End_A(S)$ is hereditary. In this case, $D(S_H)$ is a tilting $H$-module such that $A=\End_H(D(S))$ and $\Da$ is the slice determined by $D(S)$.

\end{enumerate}

\end{TheoR}

\medskip

\noindent{\sc Remark.} In case $A$ is connected, in view of Lemma \ref{lem-wcm}, we see that a slice $\Da$ of $\GaA$ is necessarily connected. As a consequence, $\Da$ is contained in a connected component $\mathcal{C}$ of $\GaA$, which is actually a section in $\mathcal{C}$; see \cite[(7.1)]{HaR} and \cite[(4.2)]{Rin}.
Such a connected component of $\GaA$ is called a {\it connecting component}.

\medskip

\section{Main results}

\medskip

In this section, we shall present our main results, that is to characterize tilted algebras in terms of the notion of a cut as defined below and to show how to obtain some tilted quotient algebras from a given algebra.

\medskip

\begin{Defn}\label{def-presec}

A full subquiver $\Da$ of $\GaA$ is called a {\it cut} if, for each arrow $X\to Y$, the following conditions are verified.

\vspace{-1pt}

\begin{enumerate}[(1)]

\item If $X\in \Da$, then either $Y$ or $\tau Y$, but not both,  belongs to $\Da$.

\item If $Y\in \Da$, then either $X$ or $\taum X$, but not both,  belongs to $\Da$.

\end{enumerate}

\end{Defn}

\medskip

\noindent{\sc Remark.} It is known that a section in a connected component of $\GaA$ is a cut; see \cite[(2.2)]{Liu2}. For this reason, a cut is called a {\it presection} in \cite[(1.3)]{ABS}.

\medskip

\noindent{\sc Example.} (1) If $\Ga$ is a connected component of $\GaA$ which is a non-homogeneous stable tube, then each ray or coray is a cut. It is evident that a stable tube contains no section.

(2) Let $A$ be an algebra with radical squared zero given by the quiver
$$
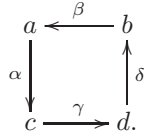
\xymatrix{
a\ar[d]_\alpha & b\ar[l]_\beta\\
c\ar[r]^\gamma & d.\ar[u]_\delta}
\vspace{5pt}$$
Its Auslander-Reiten quiver $\GaA$ is as follows:\vspace{1pt}
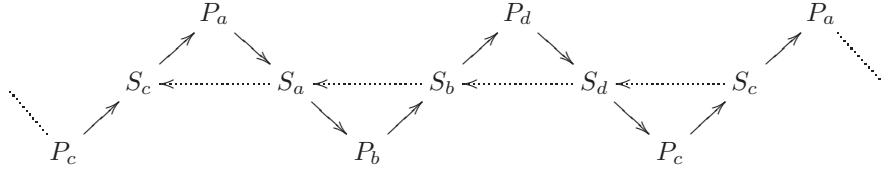
$$\hspace{-0pt}\xymatrixrowsep{12pt}\xymatrixcolsep{12pt}\xymatrix{
&&&&P_a\ar[dr] &&&& P_d \ar[dr] &&&& P_a\ar@{.}[dr] \\
&\ar@{.}[dr]&&S_c \ar[ur] && S_a\ar[dr]\ar@{.>}[ll] && S_b\ar@{.>}[ll]\ar[ur] && S_d \ar[dr]\ar@{.>}[ll] &&S_c\ar[ur]\ar@{.>}[ll] &&\\
&&P_c\ar[ur]&&&& P_b\ar[ur]&&&& P_c\ar[ur]
} \vspace{3pt}$$
where $P_a, P_b, P_c, P_d$ are the indecomposable projective modules, and $S_a, S_b, S_c, S_d$ are the simple modules. We see that $\xymatrix{\Da: \, P_b\ar[r]& S_b\ar[r]& P_d}$ is a cut in $\GaA$. Observe that $\Da$ does not meet the $\tau$-orbits of $P_a$ and $P_c$. 

\medskip

The following result exhibits the relation between a slice and a cut.

\medskip

\begin{Lemma}\label{cut-slice}

A full subquiver of $\GaA$ is a slice if and only if it is a cut which is sincere and convex in $\ind A$.

\end{Lemma}

\noindent{\it Proof.} Let $\Da$ be a full subquiver of $\GaA$. Suppose that $\Da$ is a sincere cut which is convex in $\ind A$.
In particular, $\Da$ satisfies the condition stated in Definition \ref{def-slice}(3). Let $X$ be a module in $\Da$. If $\tau X$ also belongs to $\Da$, then $\GaA$ admits a path $\tau X\to Y\to X$. By the convexity of $\Da$, we have $Y\in \Da$, a contradiction to Definition \ref{def-presec}(1). Therefore, $\Da$ is a slice.

Conversely, suppose that $\Da$ is a slice. In particular, $\Da$ satisfies the condition stated in Definition \ref{def-presec}(2). Let $X\to Y$ be an arrow in $\GaA$ with $X\in \Da$. If $Y$ is not projective, then $\tau Y \to X$ is an arrow in $\GaA$, and by Definition \ref{def-slice}(3), either $\tau Y$ or $Y$ belongs to $\Da$. Otherwise, since $\Da$ is sincere, there exists a non-zero map $f: Y\to M$ for some $M\in \Da$. By the convexity of $\Da$ in $\ind A$, we have $Y\in \Da$. This shows that $\Da$ alos satisfies the condition stated in \ref{def-presec}(1). The proof of the lemma is completed.

\medskip

Let $\Ga$ be a connected component of $\GaA$. We say that $\Ga$ is {\it semi-regular} if it contains no projective module or no injective module. In case $\Ga$ contains no oriented cycle, one says that $\Ga$ is {\it preprojective} (respectively, {\it preinjective}) if every $\tau$-orbit in $\Ga$ contains a projective (respectively, injective) module.
The following result seems unknown in the existent literature.

\medskip

\begin{Prop}\label{sr-conn}

Let $A$ be an artin algebra, and let $\mathcal{C}$ be a sincere preprojective or preinjective component of $\GaA$. If $\mathcal{C}$ is semi-regular, then $A$ is tilted with $\mathcal{C}$ being a connecting component.

\end{Prop}

\noindent{\it Proof.} We shall only consider the case where $\mathcal{C}$ is preinjective without projective modules. Then, $\mathcal{C}$ contains a section $\Da$; see \cite[(2.4)]{Liu1}. In particular, $\Da$ is a finite cut of $\GaA$.
Let $f: X\to Y$ be a non-zero map with $X\in \mathcal{C}$ and $Y\in \GaA$. Since $\mathcal{C}$ has no oriented cycle and $X$ has only finitely many successors in $\mathcal{C}$, we deduce from Lemma \ref{inf-rad}(2) that ${\rm dp}(f)<\infty.$ By Lemma \ref{fdp-map}, $\ind A$ has a non-zero path $X\rightsquigarrow Y$ of irreducible maps, and hence, $Y\in \mathcal{C}$.
Making use of this fact and the convexity of $\Da$ in $\mathcal{C}$, we see that $\Da$ is convex in $\ind A$.

Now, let $I$ be an injective module in $\GaA$. We claim that $\Hom_A(M, I)\ne 0$ for some $M\in \Da$. Indeed, we may assume that $I\not\in \Da$. Since $\mathcal{C}$ is sincere, there exists a non-zero map $g: X\to I$ with $X\in \mathcal{C}$. As shown above, $\ind A$ has a path of irreducible maps
$$\xymatrixcolsep{20pt}\xymatrix{X=X_0\ar[r]^-{f_1}& X_1\ar[r]& \cdots \ar[r] & X_{r-1} \ar[r]^-{f_r} & X_r=I}$$ with $X_i\in \mathcal{C}$ such that $f_r\cdots f_1\ne 0$. Since $\Da$ is a section, $X_i=\tau^{n_i}M_i$, where $M_i\in \Da$ and $n_i\in \Z$ with $1\ge n_i-n_{i+1}\ge 0$; see \cite[(2.3)]{Liu1}. Since $X_r$ is injective and not in $\Da$, we have $n_r< 0$. If $n_0>0$, then $n_s=0$ for some $1\le s\le r$, and our claim follows.
Suppose that $n_0<0$. Since $\mathcal{C}$ has no projective module, every minimal right almost split map for an module in $\mathcal{C}$ is surjective. Using this, we obtain an infinite path of irreducible maps \vspace{-0pt}
$$\xymatrixcolsep{20pt}\xymatrix{\cdots \ar[r] & Y_t \ar[r]^-{h_t}& Y_{t-1}\ar[r]& \cdots \ar[r] & Y_1\ar[r]^-{h_1} & Y_0=X}$$
such that $gh_1\cdots h_t\ne 0$ for every $t\ge 1$. Since $\Da$ is a finite section, some of the $Y_t$ lies in $\Da$. This establishes our claim. By Lemma \ref{cut-slice}, $\Da$ is a slice of $\GaA$. The proof of the proposition is completed.

\medskip

\begin{Lemma}\label{lem-wsm-1}

Let $\Da$ be a finite cut of $\GaA$, which is weakly convex in $\ind A$.

\vspace{-1pt}

\begin{enumerate}[$(1)$]

\item If $X\in \Da$, then neither $\tau X$ nor $\taum X$ belongs to $\Da$.

\vspace{1pt}

\item If $X, Y\in \Da$, then $\Hom_A(X, \tau Y)=0$ and $\Hom_A(\taum X, Y)=0$.

\end{enumerate}

\end{Lemma}

\noindent{\it Proof.} (1) Let $X\in \Da$. Suppose that $\tau X\in \Da$. Consider an almost split sequence
$$\xymatrixcolsep{25pt}\xymatrix{
0\ar[r] & \tau X \ar[r]^-{(f_1, u_1)} & Y_1\oplus Z_1 \ar[r]^-{g_1\choose v_1} & X\ar[r] & 0,
}$$
where $Y_1\in \Ga_A$. By the condition stated in Definition \ref{def-presec}(1), $Y_1\not\in \Da.$ Since $\Da$ is weakly convex, $g_1 f_1=0$, and consequently, $Z_1=0$; see the corollary of \cite[(1.3)]{Liu}. Using again Definition \ref{def-presec}(1), we see that $\tau Y_1\in \Da$. This yields a sectional path $Y_1\to X$ in $\GaA$ and an irreducible monomorphism $f_1: \tau X \to Y_1$ in $\mmod A$, where $\tau X\in \Da$ and $Y_1\not\in \Da$. Assume that there exists a sectional path
$$\rho_n: \; \xymatrix{
Y_n\ar[r] & Y_{n-1} \ar[r] & \cdots \ar[r] & Y_1\ar[r] & Y_0=X}$$ in $\GaA$ and an irreducible monomorphism $f_n: \tau Y_{n-1}\to Y_n$ in $\mmod A$, where $\tau Y_{n-1}, \ldots, \tau Y_0\in \Da$, while $Y_n, \ldots, Y_1\not\in \Da$. By Definition \ref{def-presec}(1), $\tau Y_n\in \Da$. Since $f_n$ is an irreducible monomorphism, there exists an almost split sequence
$$\xymatrixcolsep{35pt}
\xymatrix{
0\ar[r] & \tau Y_n \ar[r]^-{\tiny\left(\hspace{-5pt}\begin{array}{c} f_{n+1}\\ u\\ h\end{array}\hspace{-5pt}\right)} & Y_{n+1}\oplus Z_{n+1} \oplus \tau Y_{n-1} \ar[r]^-{(g, v, f_n)} & Y_n\ar[r] & 0,
}$$
where $Y_{n+1}\in \GaA$ and $f_{n+1}: \tau Y_n\to Y_{n+1}$ is a monomorphism. Observe that \vspace{2pt}
$$\rho_{n+1}: \; \xymatrix{Y_{n+1}\ar[r] & Y_n\ar[r] & Y_{n-1}\ar[r] & \cdots \ar[r] & Y_1\ar[r] & Y_0=X}\vspace{1pt}$$ is a presectional path in $\GaA$. Since $Y_n\not\in\Da$ and $Y_0\in \Da$, we deduce from Lemma \ref{pre-sec-path}(2) that $Y_{n+1}\not\in \Da$. In particular, $Y_{n+1}\ne \tau Y_{n-1}$, and thus, $\rho_{n+1}$ is a sectional path in $\GaA$. By induction, we obtain an infinite sectional path

$$\xymatrix{
\cdots \ar[r] &Y_n\ar[r] & Y_{n-1} \ar[r] & \cdots \ar[r] & Y_1\ar[r] & Y_0=X}$$ in $\GaA$ such that $\tau Y_i\in \Da,$ for all $i\ge 0.$ This yields an infinite sectional path
$$\xymatrix{
\cdots \ar[r] &\tau Y_n\ar[r] & \tau Y_{n-1} \ar[r] & \cdots \ar[r] & \tau Y_1\ar[r] & \tau Y_0}$$ in $\Da$,
a contradiction to Lemma \ref{lem-psp}. Thus, $\tau X\not\in \Da$, and consequently, $\taum X\not\in \Da$.

(2) We shall prove only the first part of Statement (2). Suppose on the contrary that there exists a non-zero map $f_0: X\to \tau Y_0,$ where $X, Y_0\in \Da$. By Statement (1), $\tau Y_0\not\in \Da$. Consider an almost split sequence
$$\xymatrix{0\ar[r] & \tau Y_0 \ar[r]^g & Z\ar[r]^h & Y_0\ar[r] & 0.}$$ Since $g$ is a monomorphism, we may find an irreducible map $g_1: \tau Y_0 \to Z_1$ with $Z_1\in \GaA$ such that $g_1f_0\ne 0$.
Since $\Da$ is weakly convex, $Z_1\not\in \Da$. Then, by Definition \ref{def-presec}(2), $Y_1=\taum Z_1\in \Da$. Continuing this process, we obtain an infinite path of irreducible maps
$$\xymatrix{\tau Y_0\ar[r]^{g_1} & \tau Y_1 \ar[r] & \cdots \ar[r] & \tau Y_{n-1} \ar[r]^{g_n} & \tau Y_n \ar[r] & \cdots, }$$ with $g_n\cdots g_1f_0\ne 0,$ where $Y_i\in \Da$ and $\tau Y_i\not\in \Da,$ for every $i\ge 0.$ Since $\Da$ is finite, we have a contradiction to Lemma \ref{lem-psp}. The proof of the lemma is completed.

\medskip

Using the following result, one can easily check whether a finite cut of $\GaA$ is weakly convex in $\ind A$.

\medskip

\begin{Prop}\label{presec-wsm}

Let $\Da$ be a cut of $\GaA$. The following conditions are equivalent.

\vspace{-2pt}

\begin{enumerate}[$(1)$]

\item The subquiver $\Da$ is finite and weakly convex in $\ind A.$

\vspace{1pt}

\item $\Hom_A(X, \tau Y)=0,$ for all $X, Y\in \Da.$

\vspace{1pt}

\item $\Hom_A(\taum\hspace{-1pt}X, Y)=0,$ for all $X, Y\in \Da.$

\end{enumerate}

\vspace{-2pt}

\noindent In this case, moreover, $\Da$ contains no oriented cycle.

\end{Prop}

\noindent{\it Proof.} Suppose first that $\Da$ is finite and weakly convex in $\ind A.$ By Lemma \ref{lem-wsm-1}, Statements (2) and (3) hold. Furthermore, assume that $\Da$ has an oriented cycle
$$\xymatrixcolsep{22pt}\xymatrix{X_0\ar[r] & \cdots \ar[r] & X_{s-1}\ar[r] &X_s=X_0.}$$
Setting $X_{s+1}=X_1$, we have $X_{i-1}=\tau X_{i+1}$ for some $1\le i\le s+1;$ see \cite{BaS}. That is, $X_{i+1}, \;\tau X_{i+1}\in \Da$, a contradiction to Lemma \ref{lem-wsm-1}(1). 

Conversely, assume that $\Hom_A(X, \tau Y)=0,$ for all $X, Y\in \Da$. It is well known that $\Da$ is finite; see \cite[Lemma 2]{Sko1}.
Suppose that there exist $L, N\in \Da$ such that $\rad^\infty(L, N)\ne 0.$ Given any $n>0$, by Lemma \ref{inf-rad}(1), $\ind A$ has a non-zero path
$$\xymatrixcolsep{22pt}\xymatrix{L\ar[r]^{u_n}& N_n\ar[r]^-{f_n} & N_{n-1}\ar[r] & \cdots \ar[r] & N_1 \ar[r]^{f_1} & N,}$$ where $N_1, \ldots, N_n\in \GaA$ and $f_1, \ldots, f_n$ are irreducible. If $N_1\not\in \Da$ then, by Definition \ref{def-presec}(2), $Y=\tau^-\hspace{-1.5pt}N_1\in \Da$, which is absurd since $\Hom_A(L, \tau Y)=0$. Thus, $N_1\in \Da$. By induction, $N_1, \ldots, N_n\in \Da$, a contradiction to Lemma \ref{lem-psp}. This shows that $\rad^\infty(X, Y)=0$, for all $X, Y\in \Da.$ Now, let
$$\xymatrixcolsep{22pt}\xymatrix{X_0\ar[r]^{g_1}& X_1 \ar[r] & \cdots \ar[r] & X_{s-1}\ar[r]^{g_s} &X_s}$$
be a non-zero path $\ind A$, where $X_0, X_s\in \Da$. Since $\rad^\infty(X_0, X_s)=0$, we have ${\rm dp}(g_1\cdots g_s)<\infty$, and consequently, ${\rm dp}(g_i)<\infty$, for $i=1, \ldots, s.$ Applying Lemma \ref{fdp-map}, we obtain a non-zero path of irreducible maps
$$\xymatrixcolsep{22pt}\xymatrix{X_0=Y_0\ar[r]^-{h_1}& Y_1 \ar[r] & \cdots \ar[r] & Y_{t-1}\ar[r]^-{h_t} &Y_t=X_s}$$ with $\{X_1, \ldots, X_{s-1}\}\subseteq \{Y_1, \ldots, Y_{t-1}\}\subseteq \GaA.$ If $Y_{t-1}\not \in \Da$, then $Z=\tau^-Y_{t-1}\in \Da.$
This yields $0\ne h_1\cdots h_{t-1}\in \Hom_A(Y_0, \tau Z)$, a contradiction. Therefore, $Y_{t-1}\in \Da$. By induction, $Y_1, \ldots, Y_{t-1}\in \Da$. In particular, $X_1, \ldots, X_{s-1}\in \Da$.
That is, $\Da$ is weakly convex in $\ind A$. Similarly, we may show that Statement (3) implies Statement (1). The proof of the proposition is completed.

\medskip

Now, we are ready to state our main result, which generalizes the result stated in \cite[(3.7)]{Ass}, \cite[(1.6)]{Liu3}, \cite[(4.2)]{Rin} and \cite[Theorem 3]{Sko2}. The proof is a refinement of the argument given in \cite[(1.6)]{Liu3}.

\medskip

\begin{Theo}\label{th-ta}

Let $A$ be an artin algebra. Then $A$ is tilted if and only if  $\GaA$ contains a faithful cut $\Da$ such that $\Hom_A(X, \tau Y)=0$ for all $X, Y\in \Da;$ and in this case, $\Da$ is a slice in $\GaA$.

\end{Theo}

{\it Proof.} If $A$ is tilted then, by Theorem \ref{Lem-slice}(2), $\GaA$ contains a finite slice $\Da$ such that the direct of its modules is a tilting module. Since tilting modules are faithful, $\Da$ is faithful. By Lemma \ref{cut-slice}, $\Da$ is a cut which is convex in $\mmod A$, and by Proposition \ref{presec-wsm}, $\Hom_A(X, Y)=0$ for all $X, Y\in \Da$.

Conversely, let $\Da$ be a faithful cut of $\GaA$ such that $\Hom_A(X, Y)=0$ for all $X, Y\in \Da$. By Proposition \ref{presec-wsm}, the direct sum $T$ of the modules in $\Da$ is a faithful module in $\mmod A$ such that $\Hom_A(T, \tau T)=0$ and $\Hom_A(\taum T, T)=0$. In particular, $\Ext_A^1(T, T)=0$; see, for example, \cite[(4.6)]{ARS}. Moreover, since $T$ is faithful, ${\rm pdim}(T)\le 1$
; see \cite[(1.5)]{RSS}.

Let $X$ be a module in $\ind A$, but not in $\Da$, such that $\Hom_A(T, X)\ne 0$. Suppose that $\Hom_A(\tau^- T, X)=0.$
Choose a non-zero map $f_0: T_0\to X$ with $T_0\in \Da$. Not being an isomorphism, $f_0$ factorizes through a minimal left almost split map $g: T_0\to L$. Therefore, there exists an irreducible map $g_1: T_0\to T_1$ with $T_1\in \GaA$ and a map $f_1: T_1\to X$ such that $f_1 g_1\ne 0.$ Since $\Hom_A(\tau^-T, X)=0,$ we have $\tau T_1\not\in \Da$, and hence, $T_1\in \Da$. By induction, we can find an infinite path of irreducible maps
$$\xymatrix{T_0\ar[r]^{g_1} & T_1 \ar[r] & \cdots \ar[r] & T_{n-1} \ar[r]^{g_n} & T_n \ar[r] & \cdots}$$
with $T_n\in \Da$ and maps $f_n: T_n\to X$ such that $g_1\cdots g_n f_n\ne 0$, for every $n\ge 1$. This is contrary to Lemma \ref{lem-psp}. Therefore, $\Hom_A(\tau^-T, X)\ne 0.$ As a consequence, $T$ is a tilting module; see \cite[(1.6)]{RSS}.

Let $H=\End_A(T)$. We claim that $H$ is hereditary. Indeed, $T$ determines a torsion theory $(\mathscr{F}, \mathscr{T})$ in $\mmod A$ with torsion class $\mathscr{T}$, and a torsion theory $(\mathscr{Y}, \mathscr{X})$ in $\mmod H$ with torsion-free class $\mathscr{Y}$.
By the Butler-Brenner Theorem; see \cite{BrB}, $\Hom_A(T, -)$ induces an equivalence from $\mathscr{T}$ to $\mathscr{Y}$. Let $P$ be an indecomposable projective module in $\mmod H$ and $v: U\to P$ be a monomorphism. Then, there exists a map $u: Q\to U$ in $\mmod H$, where $Q$ is indecomposable and projective, such that $vu\ne 0$. Observing that $u, v$ are in $\mathscr{Y}$, we may assume that $\mathscr{T}$ has morphisms $f: M\to Z$ et $h: Z\to N$, where $M, N\in \Da$ and $Z\in \ind A$, such that $Q=\Hom_A(T, M)$, $U=\Hom_A(T, Z),$ $P=\Hom_A(T, N),$ and $vu=\Hom_A(T, hf)$. Since $vu\ne 0$, we have $hf\ne 0$. By Proposition \ref{presec-wsm}, $\Da$ is weakly convex in $\ind A$. Thus, $Z\in \Da$, and hence, $U$ is projective. This establishes our claim. By Theorem \ref{Lem-slice}(2), $A$ is a tilted algebra and $\Da$ is the slice determined by $D(T_H)$. The proof of the theorem is completed.

\medskip

\noindent {\sc Remark.} As shown by the example below, the faithfulness of the cut in Theorem \ref{th-ta} cannot be replaced by the sincereness.

\medskip

\noindent{\sc Example.} Let $A$ be an algebra with radical squared zero given by the quiver
$$
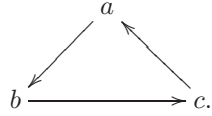
\xymatrix{&a\ar[ld] & \\
b \ar[rr] && c.\ar[ul]}
\vspace{5pt}$$
Its Auslander-Reiten quiver $\GaA$ is as follows: \vspace{-3pt}

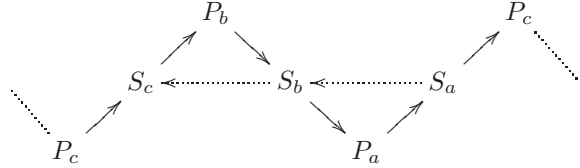
$$\hspace{-0pt}\xymatrixrowsep{12pt}\xymatrixcolsep{12pt}\xymatrix{
&&&&P_b\ar[dr] &&&& P_c \ar@{.}[dr] &&&& \\
&\ar@{.}[dr]&&S_c \ar[ur] && S_b\ar[dr]\ar@{.>}[ll] && S_a\ar@{.>}[ll]\ar[ur] && \\
&&P_c\ar[ur]&&&& P_a\ar[ur]&&&&
} \vspace{3pt}$$
where $P_a, P_b, P_c$ are the indecomposable projective modules, and $S_a, S_b, S_c$ are the simple modules. It is easy to see that $\xymatrix{\Da: P_b\ar[r]& S_b\ar[r]& P_a}$ is a sincere cut in $\GaA$ such that $\Hom_A(X, \tau Y)=0$ for all $X, Y\in \Da$. However, $A$ is not tilted.

\medskip

Applying Theorem \ref{th-ta}, we obtain the following result, which strictly includes the corresponding results stated in \cite[(3.5)]{Ass}, \cite[(2.2)]{Liu3}, \cite[(2.7)]{Liu4} and \cite[(3.1),(3.2)]{Sko3}.

\medskip

\begin{Theo}\label{tm-qta}

Let $A$ be an artin algebra, and let $\Da$ be a cut of $\Ga\hspace{-2pt}_A$ such that $\Hom_A(X, \tau Y)=0$ for all $X, Y\in \Da.$ Then $B\hspace{-2pt}=\hspace{-1pt}A/{\rm ann}(\Da)$ is a tilted algebra with $\Da$ being a slice of $\Ga_B$.

\end{Theo}

\noindent{\it Proof.} We shall identify $\mmod B$ with the full subcategory of $\mmod A$ generated by the modules annihilated by the ideal ${\rm ann}(\Da)$. In particular, $\Da$ is a faithful full subquiver of $\Ga_B$. Given any minimal right almost split map $v: E\to X$ in $\mmod A$ with $X\in \Da,$ we claim that the following statements hold true.

\vspace{-2pt}

\begin{enumerate}

\item The map $v$ is a minimal right almost map in $\mmod B.$

\item If $X$ is projective in $\mmod A$, then it is projective in $\mmod B$. Otherwise, we have $\tau_{\hspace{-1pt}_B}X=\tau X$.

\end{enumerate}

Indeed, suppose first that $v: E\to X$ is a minimal right almost map in $\mmod B.$ If $X$ is projective in $\mmod A$, then $v$ is a monomorphism, and consequently, $X$ is projective in $\mmod B$. Otherwise, there exists an almost split sequence $$\xymatrixcolsep{22pt}\xymatrix{0\ar[r] & \tau X \ar[r]^u & E \ar[r]^v & X \ar[r] & 0}$$ in $\mmod A.$ Since $u$ is a monomorphism, $\tau X\in \mmod B$, and hence, the sequence is an almost split sequence in $\mmod B$. In particular, $\tau_{_B}X=\tau X$. In order to prove Statement (1), observe that $\Da$ has no oriented cycle by Proposition \ref{presec-wsm}. Thus, the number of paths in $\Da$ starting with $X$ is finite, and we denote by $d_X$ the maximal length of such paths. If $d_X=0$, then $X$ is a sink in $\Da$. By Definition \ref{def-presec}(2), every arrow $Y\to X$ in $\GaA$ lies in $\Da$, and hence, in $\Ga_B$. As a consequence, the map $v$ lies in $\mmod B$, and hence, it is a minimal right almost split map in $\mmod B$. Assume that $d_X>0$. Let $Y\to X$ be an arrow in $\GaA$. If $Y\in \Da$, then $Y\to X$ is an arrow in $\Ga_B$. Otherwise, $Z=\tau^-Y\in \Da$ with $d_Z<d_X$. By the induction hypothesis, $Y=\tau_{\hspace{-1pt}_B}Z$ is a $B$-module, and thus, $Y\to X$ is an arrow in $\Ga_B$. Therefore, $v: E\to X$ is a minimal right almost map in $\mmod B.$ Our claim is established.

Now, by Statement (2), $\Hom_B(Y, \tau_{\hspace{-1pt}_B}X)=0$ for all $Y, X\in \Da$. Moreover, let $M\to N$ be an arrow in $\Ga_B$ with $N\in \Da$ and $M\not\in \Da$. By Statement (1), $M\to N$ is an arrow in $\GaA$. Since $\Da$ is a cut of $\GaA$, we have $\taum M\in \Da$. Using again Statement (2), we see that $\tau_{\hspace{-1pt}_B}^-M=\taum M$. This shows that $\Da$, as a subquiver of $\Ga_B$, satisfies the condition stated in Definition \ref{def-presec}(2), and dually, it also satisfies the condition stated in Definition \ref{def-presec}(1). That is, $\Da$ is a cut of $\Ga_B$. By Theorem \ref{th-ta}, $B$ is a tilted algebra with $\Da$ being a slice of $\Ga_B$. The proof of the theorem is completed.

\medskip

\noindent{\sc Remark.} Let $A$ be a cluster-tilted algebra. If $\Sa$ is a local slice in $\GaA$, by Theorem 19 stated in \cite{ABS}, $A/{\rm ann}(\Sa)$ is a tilted algebra. This fact can be deduced from Theorem \ref{tm-qta}. Indeed, in this situation, $\Sa$ is a cut such that $\Hom_A(X, \tau Y)=0$ for all $X, Y\in \Sa$; see \cite[Lemma 9]{ABS}.

\medskip

\noindent{\sc Example.} We consider again the algebra $A$ with radical squared zero given by the following quiver
$$
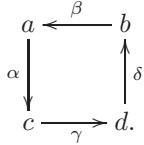
\xymatrix{
a\ar[d]_\alpha & b\ar[l]_\beta\\
c\ar[r]_\gamma & d.\ar[u]_\delta}
\vspace{5pt}$$
The Auslander-Reiten quiver $\GaA$ is as follows:
\vspace{-5pt}

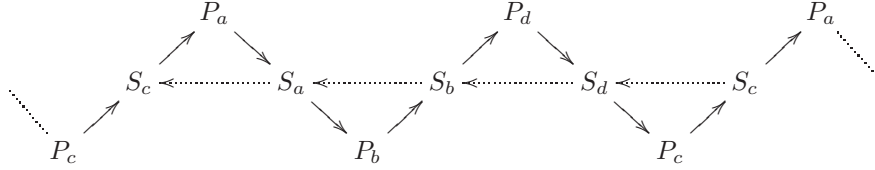
$$\hspace{-0pt}\xymatrixrowsep{12pt}\xymatrixcolsep{12pt}\xymatrix{
&&&&P_a\ar[dr] &&&& P_d \ar[dr] &&&& P_a\ar@{.}[dr] \\
&\ar@{.}[dr]&&S_c \ar[ur] && S_a\ar[dr]\ar@{.>}[ll] && S_b\ar@{.>}[ll]\ar[ur] && S_d \ar[dr]\ar@{.>}[ll] &&S_c\ar[ur]\ar@{.>}[ll] &&\\
&&P_c\ar[ur]&&&& P_b\ar[ur]&&&& P_c\ar[ur]
} \vspace{3pt}$$
We have seen that $\xymatrix{\Da: \, P_b\ar[r]& S_b\ar[r]& P_d}$ is a cut in $\GaA$. Now, it is easy to verify that $\Hom_A(X, \tau Y)=0$ for all $X, Y\in \Da$. Observe that ${\rm ann}(\Da)=Ae_cA$ and $B=A/{\rm ann}(\Da)$ is given by the  quiver

\vspace{-5pt}

$$\xymatrix{ d\ar[r]^\delta & b\ar[r]^\beta & a.} \vspace{5pt}$$
with relation $\beta\delta$. Clearly, $B$ is tilted of type $\mathbb{A}_3$. Note that $\Da$ is neither a left section nor a local slice. Therefore, none of the corresponding results stated in \cite[(3.5)]{Ass}, \cite[(2.2)]{Liu3}, \cite[(2.7)]{Liu4} and \cite[(3.1),(3.2)]{Sko3} is applicable in this example.

\end{document}